\documentclass[11pt]{article}
\usepackage[cp1251]{inputenc}
\usepackage[english]{babel}
\usepackage{amsmath}
\usepackage{amsfonts}
\usepackage{amssymb}
\usepackage{graphicx}
\newtheorem{theorem}{Theorem}
\newtheorem{proposition}{Proposition}
\newtheorem{lemma}{Lemma}

\newenvironment{proof}[1][Proof. ]{\textbf{#1}}{$\square$}
\def\text{\hbox} 
\begin{document}

\title{
\date{}
{
\begin{flushleft}
\end{flushleft}
\large \textbf{Volume formula for a $\mathbb{Z}_2$-symmetric spherical tetrahedron through its edge lengths}
\thanks{supported by the Swiss National Science Foundation no.~200020-113199/1 and RFBR no.~09-01-00255, RFBR no.~10-01-00642}}}
\author{\small Alexander Kolpakov, Alexander Mednykh, Marina Pa\v{s}kevi\v{c}}
\maketitle

\begin{abstract}\noindent
The present paper considers volume formul{\ae}, as well as trigonometric identities, that hold for a tetrahedron in 3-dimensional spherical space of constant sectional curvature $+1$. The tetrahedron possesses a certain symmetry: namely rotation through angle $\pi$ in the middle points of a certain pair of its skew edges.\\

\end{abstract}

\section{Introduction}

The volume calculation problem stated for a three-dimensional polyhedra is one of the most hard and old problems in the field of geometry. The first results belonging to the field is due to N.~Fontana~Tartaglia (1499-1557), who found a formula for the volume of a Euclidean tetrahedron better know in the present time as the Cayley-Menger determinant. Due to the paper \cite{Sabitov} the volume of every Euclidean polyhedron is a root of an algebraic equation depending on its combinatorial type and metric parameters.

In case of hyperbolic and spherical spaces the task becomes harder. The formul{\ae} for volumes
of orthoschemes are known since the work of N.~Lobachevsky and L.~Schl\"afli. The volumes of hyperbolic polyhedra with at least one ideal vertex and hyperbolic Lambert cube are given in the papers \cite{Vinberg, Kellerhals, Mednykh}. 

The general formula for the volume of a non-euclidean tetrahedron were given in \cite{ChoKim, MurakamiYano} as a linear combination of dilogarithm functions depending on the dihedral angles. Afterwards, an elementary integral formula was suggested in \cite{DM2005}. 

In case given polyhedron possesses certain symmetries, the volume formul{\ae} become facile. First, this fact was noted by Lobachevsky for ideal hyperbolic tetrahedra: the vertices of such tetrahedra belong to the ideal boundary of hyperbolic space and dihedral angles along every pair of skew edges are equal.  Later, J.~Milnor presented the respective result in very elegant form \cite{Milnor}. The general case of a tetrahedron with the same kind of symmetry is considered in \cite{DMP}. For more complicated polyhedra one also expects the use of their symmetries to be an effective tool. The volumes of octahedra enjoying certain symmetries were computed in \cite{AGM}.

The volume formula for a hyperbolic tetrahedron in terms of its edge lengths instead of dihedral angles was suggested first by \cite{MurakamiUshijima} in view of the Volume Conjecture due to R.~Kashaev \cite{Kashaev}. The further investigation on the subject was carried out in \cite{Luo2008}. 

The paper \cite{MurakamiYano} suggests a volume formula for a spherical tetrahedron as an analytic continuation of the given volume function for a hyperbolic one. The corresponding analytical strata has to be chosen in the unique proper manner. 

The present papers provides volume formula for a spherical tetrahedron that is invariant up to isometry under rotation through angle $\pi$ in the middle points of a certain pair of its skew edges. The formula itself depends on the edge lengths of given tetrahedra as well as on its dihedral angles and specifies the actual analytic strata of the volume function. Volumes of spherical Lambert cube and spherical octahedra with various kinds of symmetry were obtained in \cite{AGM, DM2009}. The analytic formul{\ae} for these polyhedra are of simpler form in contrast to their more complicated combinatorial structure.

\section{Preliminary results}

Let $\mathbb{R}^{n+1} = \{ \mathrm{x} = (x_0,\dots,x_n) : x_i \in \mathbb{R}, i=1,\dots,n \}$ be Euclidean space equipped with the standard inner product $\langle \mathrm{x}, \mathrm{y} \rangle = \sum_{i=0}^n x_i y_i$ and norm $\| \mathrm{x} \| = \sqrt{\langle \mathrm{x}, \mathrm{x} \rangle}$. Let $\mathrm{p}_0,\,\dots,\,\mathrm{p}_n$ be vectors in $\mathbb{R}^{n+1}$. Define a \textit{cone} over a collection of vectors $\mathrm{p}_0,\,\dots,\,\mathrm{p}_n$ as
\begin{equation*}
\mathrm{cone}\,\{ \mathrm{p}_0,\,\dots,\,\mathrm{p}_n \} = \left\{\sum_{i=0}^n \lambda_i \mathrm{p}_i : \lambda_i \geq 0, i=1,\dots,n\right\}.
\end{equation*}

Call a spherical $n$-simplex $\mathcal{S}$ the intersection of the cone over collection $\mathrm{p}_0,\,\dots,\,\mathrm{p}_n$ of linearly independent unitary vectors and the $n$-dimensional sphere $\mathbb{S}^n = \{v\in\mathbb{R}^n : \langle v,v \rangle = 1\}$. Thus 
\begin{equation*}
\mathcal{S} = \mathrm{cone}\,\{ \mathrm{p}_0,\,\dots,\,\mathrm{p}_n \} \cap \mathbb{S}^n.
\end{equation*}
Call $\mathrm{p}_0,\,\dots,\,\mathrm{p}_n$ the vertices of a simplex $\mathcal{S}$. Notice that $\{\mathrm{p}_0,\,\dots,\,\mathrm{p}_n\} \subset \mathbb{S}^n$. Call a $(k-1)$-dimensional face of $\mathcal{S}$ the intersection of the cone over a $k$-element sub-collection of linearly independent vectors $\{ \mathrm{p}_{i_1},\dots,\mathrm{p}_{i_k}\} \subset \{\mathrm{p}_0,\,\dots,\,\mathrm{p}_n\}$ with $\{i_1 < \dots < i_k\}\subset\{0,\dots,n\}$, $0\leq k \leq n$ and the sphere $\mathbb{S}^n$.

The matrix $G^{\star}=\{\langle \mathrm{p}_i, \mathrm{p}_j \rangle\}_{i,j=0}^n$ is the \textit{edge matrix} of a simplex $\mathcal{S}$. 

Let $M=\{ m_{ij} \}_{i,j=0}^n$ be a matrix. Denote by $M(i,j)$ the matrix obtained from $M$ by deletion of $i$-th row and $j$-th column for $i,j = 0,\dots,n$. Put $M_{ij} = (-1)^{i+j} \det M(i,j)$. The quantity $M_{ij}$ is the $(i,j)$-cofactor of $M$. Then, call $\mathrm{cof}\,M = \{ M_{ij} \}_{i,j=0}^n$ the \textit{cofactor matrix} of a matrix $M$.

The unit (outer) normal vector $\mathrm{v}_i$, $i=0,\dots,n$ to a $(n-1)$-dimensional face $\mathcal{S}_i = \{ \mathrm{p}_{i_1},\dots,\mathrm{p}_{i_n} \} \nsupseteq \{\mathrm{p}_{i}\}$ of the simplex $\mathcal{S}$ is defined by (cf. \cite{Luo1997})
\begin{equation*}
\mathrm{v}_i = \frac{\sum_{k=0, k\neq i}^n G^{\star}_{ik}\mathrm{p}_k}{\sqrt{G^{\star}_{ii}\,\det G^{\star}}}.
\end{equation*}
The matrix $G=\{\langle \mathrm{v}_i, \mathrm{v}_j \rangle\}_{i,j=0}^n$ is the \textit{Gram matrix} of a simplex $\mathcal{S}$.

Given a simplex $\mathcal{S} \subset \mathbb{S}^n$ with vertices $\{ \mathrm{p}_1,\dots,\mathrm{p}_n\}$ and unit (outer) normal vectors $\{\mathrm{v}_1,\dots,\mathrm{v}_n\}$ define its edge lengths by $\cos l_{ij} = \langle \mathrm{p}_i,\mathrm{p}_j \rangle$ and (inner) dihedral angles by  $\cos \alpha_{ij} = - \langle \mathrm{v}_i,\mathrm{v}_j \rangle$ with $0 \leq l_{ij},\,\alpha_{ij} \leq \pi$, $i,j=0,\dots,n$. Then the \textit{Gram matrix} of $\mathcal{S}$ is $G = \{-\cos \alpha_{ij} \}^n_{i,j=0}$ and the \textit{edge matrix} of $\mathcal{S}$ is $G^{\star} = \{\cos l_{ij} \}^n_{i,j=0}$.

The sphere $\mathbb{S}^n$ is endowed with the natural metric of constant sectional curvature $+1$. Call the given metric space \textit{the spherical space} $\mathbb{S}^n$. The isometry group of the spherical space $\mathbb{S}^n$ is the orthogonal group $O(n+1)$. Orientation-preserving isometries of $\mathbb{S}^n$ compose the subgroup of index two in $O(n+1)$, called $SO(n+1)$.

The following theorems tackle existence of a spherical simplex with given Gram matrix or edge matrix \cite{Luo1997}:

\begin{theorem}\label{theorem_existence_angles}
The Gram matrix $\{ -\cos \alpha_{ij} \}^n_{i,j=0}$ of a spherical $n$-simplex is symmetric, positive definite with diagonal entries equal to $1$. Conversely, every positive
definite symmetric matrix with diagonal entries equal to $1$ is the Gram matrix of a
spherical $n$-simplex that is unique up to an isometry.
\end{theorem}

\begin{theorem}\label{theorem_existence_lengths}
The edge matrix $\{ \cos l_{ij} \}^n_{i,j=0}$ of a spherical $n$-simplex is symmetric,
positive definite with diagonal entries equal to $1$. Conversely, every positive
definite symmetric matrix with diagonal entries equal to $1$ is the edge matrix of a
spherical $n$-simplex that is unique up to an isometry.
\end{theorem}

The following theorem due to Ludwig Schl\"afli relates the volume of a given simplex in the spherical space $\mathbb{S}^n$ with volumes of its apexes and dihedral angles between its faces \cite{Milnor, Schlaefli}:

\begin{theorem}[Schl\"afli formula]\label{Schlaefli_formula}
Let a simplex $\mathcal{S}$ in the space $\mathbb{S}^n$, $n \geq 2$, of constant sectional curvature $+1$ have dihedral angles $\alpha_{ij} = \angle\, \mathcal{S}_i\, \mathcal{S}_j$, $0 \leq i < j \leq n$ formed by the $(n-1)$-dimensional faces $\mathcal{S}_i$ and $\mathcal{S}_j$ of $\mathcal{S}$ which intersect in the $(n-2)$-dimensional apex $\mathcal{S}_{ij} = \mathcal{S}_i \cap \mathcal{S}_j$.

Then the differential of the volume function $\mathrm{Vol}_n$ on the set of all simplices in $\mathbb{S}^n$ satisfies the equality
\begin{equation*}
(n-1)\, \mathrm{d}\mathrm{Vol}_n\,\mathcal{S} = \sum^{n}_{i<j=0} \mathrm{Vol}_{n-2}\,\mathcal{S}_{ij}\, \mathrm{d}\alpha_{ij}
\end{equation*}
where $\mathrm{Vol}_{n-1}\,\mathcal{S}_{ij}$ is the $(n-2)$-dimensional volume function on the set of all $(n-2)$-dimensional apexes $\mathcal{S}_{ij}$, $\mathrm{Vol}_{0}\,\mathcal{S}_{ij} = 1$, $0 \leq i < j \leq n$, and $\alpha_{ij}$ is the dihedral angle between $\mathcal{S}_i$ and $\mathcal{S}_j$ along $\mathcal{S}_{ij}$. 
\end{theorem}

The Schl\"afli formula for the spherical space $\mathbb{S}^3$ can be reduced to
\begin{equation*}
\mathrm{d}\mathrm{Vol}\,\mathcal{S} = \frac{1}{2} \sum^{3}_{i<j=0} l_{ij}\, \mathrm{d}\alpha_{ij},
\end{equation*}    
where $\mathrm{Vol} = \mathrm{Vol}_3$ is the volume function, $l_{ij}$ represents the length of $ij$--th edge and $\alpha_{ij}$ represents the dihedral angle along it. So the volume of a simplex in $\mathbb{S}^3$ is related with its edge lengths and dihedral angles.

Given a simplex $\mathcal{S} \subset \mathbb{S}^n$ with vertices $\{ \mathrm{p}_1,\dots,\mathrm{p}_n\}$ and unit normal vectors $\{\mathrm{v}_1,\dots,\mathrm{v}_n\}$ denote \textit{its dual} $\mathcal{S}^{\star} \subset \mathbb{S}^n$ as a simplex with vertices $\{\mathrm{v}_1,\dots,\mathrm{v}_n\}$ and unit normal vectors $\{\mathrm{v}_1,\dots,\mathrm{v}_n\}$.

In case of the spherical space $\mathbb{S}^3$ every edge $\mathrm{p}_i\mathrm{p}_j$, $0 \leq i < j \leq 3$ of $\mathcal{S}$ corresponds to the edge $\mathrm{v}_{3-j}\mathrm{v}_{3-i}$ of its dual $\mathcal{S}^{\star}$. The theorem below was originally discovered by Italian mathematician Duke Gaetano Sforza and could be found in the book \cite{Milnor}.

\begin{theorem}\label{theorem_Sforza}
Let $\mathcal{S}$ be a simplex in the spherical space $\mathbb{S}^3$ and let $\mathcal{S}^{\star}$ be its dual. Then
\begin{equation*}
\mathrm{Vol}_3\,\mathcal{S} + \mathrm{Vol}_3\,\mathcal{S}^{\star} + \frac{1}{2}\sum_{E\subset \mathcal{S}} \mathrm{Vol}_1E\, \mathrm{Vol}_1E^{\star} = \pi^2,
\end{equation*}  
where the sum is taken over all edges $E$ of $\mathcal{S}$ and $E^{\star}$ denotes the edge of $\mathcal{S}^{\star}$ corresponding to $E$.
\end{theorem}

In what follows we call $3$-dimensional simplex \textit{a tetrahedron} for the sake of brevity.


\section{Trigonometric identities for a spherical tetrahedron}

Let $\mathbf{T}$ be a tetrahedron in the spherical space $\mathbb{S}^3$ with vertices $\mathrm{p}_0$, $\mathrm{p}_1$, $\mathrm{p}_2$, $\mathrm{p}_3$, dihedral angles $A$, $B$, $C$, $D$, $E$, $F$ and edge lengths $l_A$, $l_B$, $l_C$, $l_D$, $l_E$, $l_F$. In the sequel we have that $0 \leq A,\,B,\,C,\,D,\,E,\,F \leq \pi$ and $0 \leq l_A,\,l_B,\,l_C,\,l_D,\,l_E,\,l_F \leq \pi$. 

\begin{figure}[ht]
\begin{center}
\includegraphics* [totalheight=6cm]{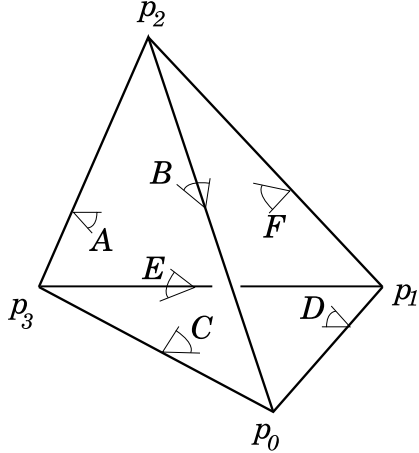}
\end{center}
\caption{Spherical tetrahedron} \label{tetr}
\end{figure}

Denote by
\begin{equation*}
G^{\star} = \{ g^{\star}_{ij} \}_{i,j=0}^3 = 
\left(\begin{array}{cccc} 1 &\cos l_A &\cos l_B &\cos l_C\\
\cos l_A &1 &\cos l_F &\cos l_E\\
\cos l_B &\cos l_F &1 &\cos l_D\\
\cos l_C &\cos l_E &\cos l_D &1 \end{array}\right)
\end{equation*}
the \textit{edge matrix} of $\mathbf{T}$ and by
\begin{equation*}
G = \{ g_{ij} \}_{i,j=0}^3 = \left(\begin{array}{cccc} 1 &-\cos D &-\cos E &-\cos F\\
-\cos D &1 &-\cos C &-\cos B\\
-\cos E &-\cos C &1 &-\cos A\\
-\cos F &-\cos B &-\cos A &1 \end{array}\right)
\end{equation*}
its \textit{Gram matrix}.

Denote by $c_{ij}$ and $c^{\star}_{ij}$ the respective cofactors of the matrices $G$ and $G^{\star}$ for $i,j = 0,1,2,3$. 

Further, we mention several important trigonometric relations (see, e.g. \cite{DMP}) to be used below:

\begin{theorem}[Sine Rule]\label{SineRule_general}
Given a spherical tetrahedron $\mathbf{T}$ with the Gram matrix $G$  and the edge matrix $G^{\star}$, denote $\Delta = \det G$, $\Delta^{\star} = \det G^{\star}$, $p = c_{00} c_{11} c_{22} c_{33}$, $p^{\star} = c^{\star}_{00} c^{\star}_{11} c^{\star}_{22} c^{\star}_{33}$. Then
\begin{equation*}
\frac{\sin l_A \sin l_D}{\sin A \sin D} = \frac{\sin l_B \sin l_E}{\sin B \sin E} = \frac{\sin l_C \sin l_F}{\sin C \sin F} = \frac{\Delta}{\sqrt{p}} = \frac{\sqrt{p^{\star}}}{\Delta^{\star}}.
\end{equation*}
\end{theorem}
\bigskip
\begin{theorem}[Cosine Rule]\label{CosineRule_general}
For the respective pairs of skew edges of a spherical tetrahedron $\mathbf{T}$ the following equalities hold:
\begin{equation*}
\frac{\cos l_A \cos l_D - \cos l_B \cos l_E}{\cos A \cos D - \cos B \cos E} = \frac{\cos l_B \cos l_E - \cos l_C \cos l_F}{\cos B \cos E - \cos C \cos F} =
\end{equation*}
\begin{equation*}
= \frac{\cos l_C \cos l_F - \cos l_A \cos l_D}{\cos C \cos F - \cos A \cos D} = \frac{\Delta}{\sqrt{p}} = \frac{\sqrt{p^{\star}}}{\Delta^{\star}}.
\end{equation*}
\end{theorem}
\bigskip

Also we need the following theorem due to Jacobi (see \cite[Th\'eor\`eme 2.5.2]{Prasolov}):
\begin{theorem}\label{theorem_Jacobi}
Let $M = \{ m_{ij} \}_{i,j=0}^n$ be a matrix, $\mathrm{cof}\,M = \{ M_{ij} \}_{i,j=0}^n$ be its cofactor matrix, $0 < k < n$ and $\sigma = \left(\begin{array}{cc} i_0 \ldots i_n \\ j_0 \ldots j_n \end{array}\right)$ an arbitrary permutation. Then
\begin{equation*}
\det \{ M_{i_p j_q} \}^{k}_{p,q=0} = (-1)^{\mathrm{sgn}\,\sigma}(\det M)^{k}\det \{ m_{i_p j_q} \}^{n}_{p,q=k}.
\end{equation*}
\end{theorem}

\section{Trigonometric identities for a $\mathbb{Z}_2$-symmetric spherical tetrahedron}

Consider a spherical tetrahedron $\mathbf{T}$ which is symmetric under rotation through angle $\pi$ about the axe that passes through the middle points of the edges $\mathrm{p}_0\mathrm{p}_1$ and $\mathrm{p}_2\mathrm{p}_3$. We call such tetrahedron $\mathbb{Z}_2$-symmetric. Note, that in this case $l_B = l_E$, $l_C = l_F$ and $B = E$, $C = F$.

\begin{figure}[ht]
\begin{center}
\includegraphics* [totalheight=6cm]{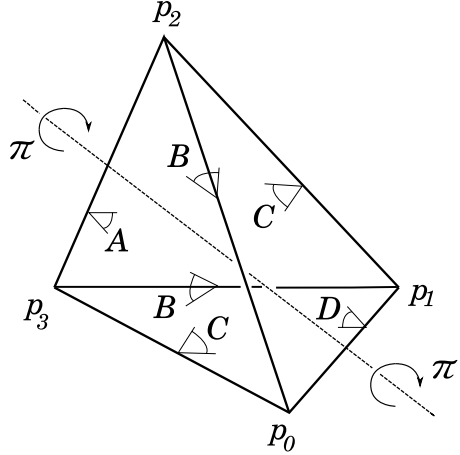}
\end{center}
\caption{$\mathbb{Z}_2$-symmetric spherical tetrahedron} \label{tetrsymm}
\end{figure}

\begin{lemma}\label{Lemma_A_D}
For a $\mathbb{Z}_2$-symmetric spherical tetrahedron with dihedral angles $A$, $B = E$, $C = F$, $D$ and edge lengths $l_A$, $l_B = l_E$, $l_C = l_F$, $l_D$ the following statements hold:
\begin{itemize}
\item[\rm{(\textbf{i})}] $l_A = l_D$ if and only if $A = D$,
\item[\rm{(\textbf{ii})}] $l_A > l_D$ if and only if $A < D$.
\end{itemize}
\end{lemma}
\begin{proof}
As applied to the edge matrix $G^{\star}$, Theorem~\ref{theorem_Jacobi} grants the equality
\begin{equation*}
c^{\star}_{00} c^{\star}_{23} - c^{\star}_{01} c^{\star}_{22} = \Delta^{\star} (g^{\star}_{01} - g^{\star}_{23}).
\end{equation*}

Recall,
\begin{equation*}
g_{01} = -\cos D = \frac{c^{\star}_{01}}{\sqrt{c^{\star}_{00} c^{\star}_{11}}},\,
g_{23} = -\cos A = \frac{c^{\star}_{23}}{\sqrt{c^{\star}_{22} c^{\star}_{33}}},
\end{equation*}
\begin{equation*}
g^{\star}_{01} = \cos l_A,\, g^{\star}_{23} = \cos l_D.
\end{equation*}

For a spherical $\mathbb{Z}_2$-symmetric tetrahedron we have $c^{\star}_{00} = c^{\star}_{11} > 0$, $c^{\star}_{22} = c^{\star}_{33} > 0$. Thus,
\begin{equation*}
\cos l_A - \cos l_D = - \frac{c^{\star}_{00} c^{\star}_{22}}{\Delta^{\star}} (\cos A - \cos D).
\end{equation*}

As far as $0 \leq l_A,\, l_D \leq \pi$ and $0 \leq A,\, D \leq \pi$, the assertions of the Lemma follow.
\end{proof}

Trigonometric identities from in Theorem~\ref{SineRule_general} and Theorem~\ref{CosineRule_general} imply 
\begin{proposition}\label{SineRule}
Let $\mathbf{T}$ be a $\mathbb{Z}_2$-symmetric spherical tetrahedron with dihedral angles $A$, $B = E$, $C = F$, $D$ and edge lengths $l_A$, $l_B = l_E$, $l_C = l_F$, $l_D$. Then the following equalities hold:
\begin{equation*}
u = \frac{\sin \frac{l_A+l_D}{2}}{\sin \frac{A+D}{2}} = \frac{\sin \frac{l_A-l_D}{2}}{\sin \frac{D-A}{2}} =
\frac{\sin l_B}{\sin B} = \frac{\sin l_C}{\sin C} = v^{-1},
\end{equation*}
where
\begin{equation*}
u = \sqrt{\frac{c^{\star}_{00} c^{\star}_{22}}{\Delta^{\star}}}, \,\,v = \sqrt{\frac{c_{00} c_{22}}{\Delta}}
\end{equation*}
correspond to the principal and the dual parameters of the tetrahedron $\mathbf{T}$.
\end{proposition}
\begin{proof}
Recall the common property of ratios:
\begin{equation*}
\frac{a}{b} = \frac{c}{d} = \frac{a-c}{b-d} = \frac{a+c}{b+d}.
\end{equation*}
By use of the equalities above and trigonometric identities
\begin{equation*}
\cos(\varphi+\psi) = \cos\varphi \sin\psi - \cos\psi \sin\varphi,
\end{equation*}
\begin{equation*}
\cos(\varphi-\psi) = \cos\varphi \sin\psi + \cos\psi \sin\varphi
\end{equation*}
we deduce from Theorem~\ref{SineRule_general} and Theorem~\ref{CosineRule_general} the following relations:
\begin{equation*}
\frac{\Delta}{c_{00} c_{22}} = \frac{1 - \cos(l_A+l_D)}{1 - \cos(A+D)} = \frac{1 - \cos(l_A-l_D)}{1 - \cos(D-A)} = \frac{\sin^2 l_B}{\sin^2 B}
= \frac{\sin^2 l_C}{\sin^2 C} = \frac{c^{\star}_{00} c^{\star}_{22}}{\Delta^{\star}}.
\end{equation*}
Put $u = \sqrt{\frac{c^{\star}_{00} c^{\star}_{22}}{\Delta^{\star}}}$, $v = \sqrt{\frac{c_{00} c_{22}}{\Delta}}$. The quantities $u$ and $v$ are positive real numbers, related together by $uv = 1$. We call $u$ \textit{the principal parameter} of the tetrahedron $\mathbf{T}$ and $v$ is called its \textit{dual parameter}.
Granted the identity $1 - \cos\varphi = 2\sin^2 \frac{\varphi}{2}$ it follows that
\begin{equation*}
u^2 = \frac{\sin^2 \frac{l_A+l_D}{2}}{\sin^2 \frac{A+D}{2}} = \frac{\sin^2 \frac{l_A-l_D}{2}}{\sin^2 \frac{D-A}{2}} = \frac{\sin^2 l_B}{\sin^2 B} = \frac{\sin^2 l_C}{\sin^2 C} = v^{-2}.
\end{equation*}
Extracting square roots in accordance with Lemma~\ref{Lemma_A_D} finishes the proof.
\end{proof}

\section{Volume formula for a $\mathbb{Z}_2$-symmetric spherical tetrahedron}

\subsection{Further trigonometric identities for a $\mathbb{Z}_2$-symmetric spherical tetrahedron}

Let $\mathbf{T}$ be a $\mathbb{Z}_2$-symmetric spherical tetrahedron with dihedral angles $A$, $B = E$, $C = F$, $D$ and edge lengths $l_A$, $l_B = l_E$, $l_C = l_F$, $l_D$. Denote
\begin{equation*}
a_{+}=\cos \frac{l_A+l_D}{2},\, a_{-}=\cos \frac{l_A-l_D}{2},\, b=\cos l_B,\, c=\cos l_C
\end{equation*}
and
\begin{equation*}
\mathcal{A}_{+} = \cos \frac{A+D}{2},\, \mathcal{A}_{-} = \cos \frac{D-A}{2},\, \mathcal{B} = \cos B,\, \mathcal{C} = \cos C.
\end{equation*}

The Lemma below gives a useful identity that follows from the definition of the principal parameter $u$ for the tetrahedron $\mathbf{T}$:
\begin{lemma}\label{lemma_u}
The principal parameter $u$ of the tetrahedron $\mathbf{T}$ is a positive root of the quadratic equation
\begin{equation*}
u^2 + \frac{4\,(a_{+}a_{-}-bc)(a_{+}b-a_{-}c)(a_{+}c-a_{-}b)}{\Delta^{\star}} = 1,
\end{equation*}
where
\begin{equation*}
\Delta^{\star} = (a_{+}+a_{-}+b+c)(a_{+}+a_{-}-b-c)(a_{+}-a_{-}-b+c)(a_{+}-a_{-}+b-c).
\end{equation*}
\end{lemma}
\begin{proof}
Substitute $u$ from Proposition \ref{SineRule}, then express the product of $c^{\star}_{00}$ and $c^{\star}_{22}$ as a polynomial in the new variables $a_{+}$, $a_{-}$, $b$, $c$ and proceed with straightforward computations.
\end{proof}

As for the dual parameter of $\mathbf{T}$, the following Lemma holds:
\begin{lemma}\label{lemma_v}
The dual parameter $u$ of the tetrahedron $\mathbf{T}$ is a positive root of the quadratic equation
\begin{equation*}
v^2 - \frac{4\,(\mathcal{A}_{+}\mathcal{A}_{-}+\mathcal{B}\mathcal{C})(\mathcal{A}_{+}\mathcal{B}+\mathcal{A}_{-}\mathcal{C})(\mathcal{A}_{+}\mathcal{C}+\mathcal{A}_{-}\mathcal{B})}{\Delta} = 1,
\end{equation*}
where
\begin{equation*}
\Delta = (\mathcal{A}_{-}-\mathcal{A}_{-}-\mathcal{B}-\mathcal{C})(\mathcal{A}_{-}-\mathcal{A}_{+}+\mathcal{B}+\mathcal{C})(\mathcal{A}_{-}+\mathcal{A}_{+}-\mathcal{B}+\mathcal{C})(\mathcal{A}_{-}+\mathcal{A}_{+}+\mathcal{B}-\mathcal{C}).
\end{equation*}
\end{lemma}

In what follows we say a spherical tetrahedron $\mathbf{T}_{s}$ with dihedral angles $\alpha$, $\beta$, $\gamma$, $\delta$, $\varepsilon$, $\varphi$  and edge lengths $l_\alpha$, $l_\beta$, $l_\gamma$, $l_\delta$, $l_\varepsilon$, $l_\varphi$ to be \textit{symmetric} if $l_\alpha = l_\delta$, $l_\beta = l_\varepsilon$, $l_\gamma = l_\varphi$ or, that is equivalent, $\alpha = \delta$, $\beta = \varepsilon$, $\delta = \varphi$.

Denote
\begin{equation*}
\tilde{a} = \cos l_\alpha,\, \tilde{b} = \cos l_\beta,\, \tilde{c} = \cos l_\gamma,
\end{equation*}
\begin{equation*}
\tilde{A} = \cos \alpha,\, \tilde{B} = \cos \beta,\, \tilde{C} = \cos \gamma,
\end{equation*}

The following trigonometric identities are proven in \cite{DMP}: 
\begin{proposition}\label{SineRule_symmetric}
Let $\mathbf{T}_s$ be a symmetric spherical tetrahedron with dihedral angles $\alpha = \delta$, $\beta = \varepsilon$, $\gamma = \varphi$  and edge lengths $l_\alpha = l_\delta$, $l_\beta = l_\varepsilon$, $l_\gamma = l_\varphi$. Then the following equalities hold:
\begin{equation*}
\frac{\sin l_\alpha}{\sin \alpha} = \frac{\sin l_\beta}{\sin \beta} = \frac{\sin l_\gamma}{\sin \gamma} = u_s,
\end{equation*}
where the positive root of quadratic equation
\begin{equation*}
u^2_s + \frac{4(\tilde{a}-\tilde{b}\tilde{c})(\tilde{b}-\tilde{a}\tilde{c})(\tilde{c}-\tilde{a}\tilde{b})}{\delta^{\star}} = 1
\end{equation*}
with
\begin{equation*}
\delta^{\star} = (\tilde{a}+\tilde{b}+\tilde{c}+1)(\tilde{a}-\tilde{b}-\tilde{c}+1)(\tilde{b}-\tilde{a}-\tilde{c}+1)(\tilde{c}-\tilde{a}-\tilde{b}+1)
\end{equation*}
represents the principal parameter $u_s$ of the tetrahedron $\mathbf{T}_s$.
\end{proposition}

Meanwhile, the dual parameter of $\mathbf{T}_s$ is the positive root of the following equation:

\begin{equation*}
v^2_s - \frac{4(\tilde{A}+\tilde{B}\tilde{C})(\tilde{B}+\tilde{A}\tilde{C})(\tilde{C}+\tilde{A}\tilde{B})}{\delta} = 1
\end{equation*}
with
\begin{equation*}
\delta = (1-\tilde{A}-\tilde{B}-\tilde{C})(1-\tilde{A}+\tilde{B}+\tilde{C})(1+\tilde{A}-\tilde{B}+\tilde{C})(1+\tilde{A}+\tilde{B}-\tilde{C}).
\end{equation*}

The following lemma shows the correspondence between $\mathbb{Z}_2$-symmetric and symmetric spherical tetrahedra.

\begin{lemma}\label{lemma_symmetric}
Let $\mathbf{T}$ be a $\mathbb{Z}_2$-symmetric spherical tetrahedron with dihedral angles $A$, $B = E$, $C = F$, $D$ and edge lengths $l_A$, $l_B = l_E$, $l_C = l_F$, $l_D$. Then there exists the associated symmetric tetrahedron $\mathbf{T}_s$ with the Gram matrix
\begin{equation*}
G_{s} = 
\left( \begin{array}{cccc}1& -\cos \alpha &-\cos \beta &-\cos \gamma\\
-\cos \alpha &1 &-\cos \gamma &-\cos \beta\\
-\cos \beta &-\cos \gamma &1 &-\cos \alpha\\
-\cos \gamma &-\cos \beta &-\cos \alpha &1 \end{array} \right)=
\end{equation*}
\begin{equation*}
=\left(\begin{array}{cccc} 1& -\frac{\mathcal{A}_{+}}{\mathcal{A}_{-}}& -\frac{\mathcal{B}}{\mathcal{A}_{-}}& -\frac{\mathcal{C}}{\mathcal{A}_{-}}\\
-\frac{\mathcal{A}_{+}}{\mathcal{A}_{-}}& 1& -\frac{\mathcal{C}}{\mathcal{A}_{-}}& -\frac{\mathcal{B}}{\mathcal{A}_{-}}\\
-\frac{\mathcal{B}}{\mathcal{A}_{-}}& -\frac{\mathcal{C}}{\mathcal{A}_{-}}& 1& -\frac{\mathcal{A}_{+}}{\mathcal{A}_{-}}\\
-\frac{\mathcal{C}}{\mathcal{A}_{-}}& -\frac{\mathcal{B}}{\mathcal{A}_{-}}& -\frac{\mathcal{A}_{+}}{\mathcal{A}_{-}}& 1 \end{array}\right)
\end{equation*}
and the edge matrix
\begin{equation*}\displaystyle
G^{\star}_{s} 
= \left( \begin{array}{cccc}1& \cos l_\alpha &\cos l_\beta &\cos l_\gamma\\
\cos l_\alpha &1 &\cos l_\gamma &\cos l_\beta\\
\cos l_\beta &\cos l_\gamma &1 &\cos l_\alpha\\
\cos l_\gamma &\cos l_\beta &\cos l_\gamma &1 \end{array} \right) 
= \left(\begin{array}{cccc} 1& \frac{a_{+}}{a_{-}}& \frac{b}{a_{-}}& \frac{c}{a_{-}}\\
\frac{a_{+}}{a_{-}}& 1& \frac{c}{a_{-}}& \frac{b}{a_{-}}\\
\frac{b}{a_{-}}& \frac{c}{a_{-}}& 1& \frac{a_{+}}{a_{-}}\\
\frac{c}{a_{-}}& \frac{b}{a_{-}}& \frac{a_{+}}{a_{-}}& 1 \end{array}\right).
\end{equation*}
\end{lemma}
\begin{proof}
Denote the principal cofactors of $G^{\star}_s$ by $s^{\star}_{ii}$, $i = 0,1,2,3$.
To prove existence of $\mathbf{T}_s$ it suffices by Theorem \ref{theorem_existence_lengths} to show that $\det G^{\star}_s > 0$ and $s^{\star}_{ii} > 0$, $i = 0,1,2,3$.
We have that $\displaystyle \det G^{\star}_s = \frac{\det G^{\star}}{a^4_{-}} > 0$ and
\begin{equation*}
s^{\star}_{00} - \frac{c^{\star}_{00}}{a^2_{-}} = - \frac{2}{a^2_{-}}(a_{+}a_{-} - bc)\sin \frac{l_A-l_D}{2} \sin l_A,
\end{equation*}
\begin{equation*}
s^{\star}_{00} - \frac{c^{\star}_{22}}{a^2_{-}} = \frac{2}{a^2_{-}}(a_{+}a_{-} - bc)\sin \frac{l_A-l_D}{2} \sin l_D.
\end{equation*}
From the former two equalities we deduce that depending on the sign of their right-hand parts either $s^{\star}_{00} \geq \frac{c^{\star}_{00}}{a^2_{-}} $ or $s^{\star}_{00} \geq \frac{c^{\star}_{22}}{a^2_{-}}$. As far as the tetrahedron $\mathbf{T}$ exists, $c^{\star}_{00} > 0$ and $c^{\star}_{22} > 0$. It follows that $s^{\star}_{00} = s^{\star}_{11} = s^{\star}_{22} = s^{\star}_{33} > 0$. Thus, the tetrahedron $\mathbf{T}_s$ exists.

For the edge lengths of the symmetric spherical tetrahedron $\mathbf{T}_s$ one has
\begin{equation*}
\cos l_\alpha = \frac{a_{+}}{a_{-}},\, \cos l_\beta = \frac{b}{a_{-}},\, \cos l_\gamma = \frac{c}{a_{-}}.
\end{equation*}
Let us prove that
\begin{equation*}
\cos \alpha = \frac{\mathcal{A}_{+}}{\mathcal{A}_{-}},\, \cos \beta = \frac{\mathcal{B}}{\mathcal{A}_{-}},\, \cos \gamma = \frac{\mathcal{C}}{\mathcal{A}_{-}}.
\end{equation*}

Refer to Proposition~\ref{SineRule} together with Lemma~\ref{lemma_u} and note that the following relation holds between the principal parameters $u$ and $u_s$ of tetrahedra $\mathbf{T}$ and $\mathbf{T}_s$, respectively:
\begin{equation*}
1 - u^2 = a_{-}^2 (1-u^2_s).
\end{equation*}
Substituting $u = \frac{\sin\frac{l_A+l_D}{2}}{\sin\frac{A+D}{2}}$ from Proposition~\ref{SineRule} and $u_s = \frac{\sin l_\alpha}{\sin \alpha}$ from Proposition~\ref{SineRule_symmetric} to the relation above one obtains
\begin{equation*}
\cos^2 \alpha = \frac{\cos^2 \frac{A+D}{2}}{1 - \frac{\sin^2\frac{l_A-l_D}{2}}{u^2}} = \frac{\cos^2 \frac{A+D}{2}}{\cos^2 \frac{D-A}{2}} = \frac{\mathcal{A}^2_{+}}{\mathcal{A}^2_{-}}.
\end{equation*}
Thus,
\begin{equation*}
\cos \alpha = \pm \frac{\cos \frac{A+D}{2}}{\cos \frac{D-A}{2}} = \pm \frac{\mathcal{A}_{+}}{\mathcal{A}_{-}}.
\end{equation*}
We should choose the proper sign in the equality above.
Note, that if $\mathbf{T}$ is symmetric, i.e. $l_A = l_D$, then $G^{\star}_s = G^{\star}$. That means $\mathbf{T}$ has isometric associated symmetric tetrahedron $\mathbf{T}_s$. Thus, the Gram matrices for tetrahedra $\mathbf{T}$ and $\mathbf{T}_s$ mentioned in the assertions of the Lemma coincide. For the equality $G_s = G$ to hold if $\mathbf{T}$ is symmetric, we put 
\begin{equation*}
\cos \alpha = \frac{\cos \frac{A+D}{2}}{\cos \frac{D-A}{2}} = \frac{\mathcal{A}_{+}}{\mathcal{A}_{-}}.
\end{equation*}
The rest of the proof follows by analogy.
\end{proof}

Denote the auxiliary parameters of the tetrahedron $\mathbf{T}$:
\begin{equation*}
t^2 = 1 - u^2 = \frac{4\,(a_{+}a_{-}-bc)(a_{+}b-a_{-}c)(a_{+}c-a_{-}b)}{\Delta^{\star}}
\end{equation*}
and
\begin{equation*}
\tau^2 = v^2 - 1 = \frac{4\,(\mathcal{A}_{+}\mathcal{A}_{-}+\mathcal{B}\mathcal{C})(\mathcal{A}_{+}\mathcal{B}+\mathcal{A}_{-}\mathcal{C})(\mathcal{A}_{+}\mathcal{C}+\mathcal{A}_{-}\mathcal{B})}{\Delta}.
\end{equation*}

By Lemma~\ref{lemma_u} the quantity $t$ could be either real or pure imaginary. We choose $t$ to be non-negative or to have non-negative imaginary part. Under the same rule the quantity $\tau$ is chosen. From Proposition~\ref{SineRule} it follows that $\tau = t/u$.

The quantity $t$ is related to the parameters $a_{+}$, $a_{-}$, $b$ and $c$ of the tetrahedron $\mathbf{T}$ in the following way:

\begin{lemma}\label{lemma_square_roots}
Let $\mathbf{T}$ be a $\mathbb{Z}_2$-symmetric spherical tetrahedron with dihedral angles $A$, $B = E$, $C = F$, $D$ and edge lengths $l_A$, $l_B = l_E$, $l_C = l_F$, $l_D$. Then
\begin{itemize}
\item[\rm{(\textbf{i})}] $a^2_{-} - t^2 = a^6_{-}\, \frac{\displaystyle (s^{\star}_{00})^2}{\displaystyle \Delta^{\star}}$,
\item[\rm{(\textbf{ii})}] $a^2_{+} - t^2 = a^6_{-}\, \frac{\displaystyle (s^{\star}_{01})^2}{\displaystyle \Delta^{\star}}$,
\item[\rm{(\textbf{iii})}] $b^2 - t^2 = a^6_{-}\, \frac{\displaystyle (s^{\star}_{02})^2}{\displaystyle \Delta^{\star}}$,
\item[\rm{(\textbf{iv})}] $c^2 - t^2 = a^6_{-}\, \frac{\displaystyle (s^{\star}_{03})^2}{\displaystyle \Delta^{\star}}$;
\end{itemize}
where $s^{\star}_{ij}$, $i, j = 0,1,2,3$ are respective cofactors of the matrix
\begin{equation*}
G^{\star}_{s} = \left(\begin{array}{cccc} 1& \frac{a_{+}}{a_{-}}& \frac{b}{a_{-}}& \frac{c}{a_{-}}\\
\frac{a_{+}}{a_{-}}& 1& \frac{c}{a_{-}}& \frac{b}{a_{-}}\\
\frac{b}{a_{-}}& \frac{c}{a_{-}}& 1& \frac{a_{+}}{a_{-}}\\
\frac{c}{a_{-}}& \frac{b}{a_{-}}& \frac{a_{+}}{a_{-}}& 1 \end{array}\right).
\end{equation*}
\end{lemma}
\begin{proof}
Substitute the expression for $t^2$ from above and proceed with straightforward computations.
\end{proof}

The following proposition is used to determine signs of the respective cofactors $c_{ij}$ and $c^{\star}_{ij}$ for $i,j=0,1,2,3$ of the matrices $G = \{ g_{ij} \}_{i,j=0}^3$ and $G^{\star} = \{ g^{\star}_{ij} \}_{i,j=0}^3$ depending on the signs of their entries:

\begin{proposition}\label{Inequality_Gram}
The following inequalities hold between entries and cofactors of Gram and edge matrices for a spherical tetrahedron $\mathbf{T}$:
\begin{itemize}
\item $g_{ij} c^{\star}_{ij} \geq 0$,
\item $g^{\star}_{ij} c_{ij} \geq 0$;
\end{itemize}
where $i,j = 0,1,2,3$.
\end{proposition}
\begin{proof}
By \cite[Ch.~1, \textsection~4.2]{Vinberg} we have
\begin{equation*}
g_{ij} = \frac{c^{\star}_{ij}}{\sqrt{c^{\star}_{ii} c^{\star}_{jj}}},\,g^{\star}_{ij} = \frac{c_{ij}}{\sqrt{c_{ii} c_{jj}}}
\end{equation*}
and $c_{ii} > 0$, $c^{\star}_{ii} > 0$ for $i, j = 0,1,2,3$.

Thus
\begin{equation*}
g_{ij} c^{\star}_{ij} = \frac{(c^{\star}_{ij})^2}{\sqrt{c^{\star}_{ii} c^{\star}_{jj}}} \geq 0,\,g^{\star}_{ij} c_{ij} = \frac{c^{2}_{ij}}{\sqrt{c_{ii} c_{jj}}} \geq 0,
\end{equation*}
where $i, j = 0,1,2,3$.
\end{proof}

The following Lemma provides some useful identities that are used below:

\begin{lemma}\label{lemma_inv_hyp_identities}
The following equalities hold:
\begin{itemize}
\item[\rm{(\textbf{i})}] $\mathrm{Re}\,\sinh^{-1}(x) + \mathrm{Re}\,\sinh^{-1}(y) = \mathrm{Re}\,\sinh^{-1}(x\sqrt{y^2-1}+y\sqrt{x^2-1})$, where $x,y \in i \mathbb{R}$, $\mathrm{Im}\,x, \mathrm{Im}\,y \geq 0$,
\item[\rm{(\textbf{ii})}] $\mathrm{Re}\,\sinh^{-1}(x) - \mathrm{Re}\,\sinh^{-1}(y) = \mathrm{Re}\,\sinh^{-1}(-x\sqrt{y^2-1}+y\sqrt{x^2-1})$, where $x,y \in i \mathbb{R}$, $\mathrm{Im}\,x, \mathrm{Im}\,y \geq 0$,
\item[\rm{(\textbf{iii})}] $\mathrm{Re}\,\sinh^{-1}(x) + \mathrm{Re}\,\sinh^{-1}(y) = \mathrm{Re}\,\sinh^{-1}(x\sqrt{y^2+1}+y\sqrt{x^2+1})$, where $x,y \in \mathbb{R}$, $x,y \geq 0$,
\item[\rm{(\textbf{iv})}] $\mathrm{Re}\,\sinh^{-1}(x) - \mathrm{Re}\,\sinh^{-1}(y) = \mathrm{Re}\,\sinh^{-1}(x\sqrt{y^2+1}-y\sqrt{x^2+1})$, where $x,y \in \mathbb{R}$, $x,y \geq 0$,
\end{itemize}
\end{lemma}
\begin{proof}
Using the logarithmic representation for the function $\sinh^{-1}(\circ)$ and properties of the complex logarithm $\log(\circ)$ one derives the statement of the Lemma for the real parts of corresponding expressions.
\end{proof}

We need the relations below to derive the volume formul{\ae} for a $\mathbb{Z}_2$-symmetric spherical tetrahedron.

\begin{proposition}\label{proposition_real_parts}
Let $\mathbf{T}$ be a $\mathbb{Z}_2$-symmetric spherical tetrahedron with dihedral angles $A$, $B = E$, $C = F$, $D$ and edge lengths $l_A$, $l_B = l_E$, $l_C = l_F$, $l_D$. Without loss of generality, assume that $l_A \geq l_D$ or, equivalently, $D \geq A$ and, furthermore, $B \leq C$. 

Then the following cases are possible:
\begin{itemize}
\item[\rm{(\textbf{i})}] if $A+D \geq \pi$, $B \geq \frac{\pi}{2}$, $C \geq \frac{\pi}{2}$ and $t^2 \leq 0$, then
\begin{equation*}
\mathrm{Re}\,\left(\sinh^{-1}\frac{a_{+}}{t}+\sinh^{-1}\frac{b}{t}+\sinh^{-1}\frac{c}{t}+\sinh^{-1}\frac{a_{-}}{t}\right) = 0.
\end{equation*}

\item[\rm{(\textbf{ii})}] if $A+D \geq \pi$, $B \leq \frac{\pi}{2}$, $C \geq \frac{\pi}{2}$, then $t^2 \geq 0$ and
\begin{equation*}
\mathrm{Re}\,\left(-\sinh^{-1}\frac{a_{+}}{t}+\sinh^{-1}\frac{b}{t}-\sinh^{-1}\frac{c}{t}-\sinh^{-1}\frac{a_{-}}{t}\right) = 0,
\end{equation*}

\item[\rm{(\textbf{iii})}] if $A+D \geq \pi$, $B \leq \frac{\pi}{2}$, $C \leq \frac{\pi}{2}$ and $t^2 \leq 0$, then
\begin{equation*}
\mathrm{Re}\, \left(-\sinh^{-1}\frac{a_{+}}{t}+\sinh^{-1}\frac{b}{t}+\sinh^{-1}\frac{c}{t}-\sinh^{-1}\frac{a_{-}}{t}\right) = 0,
\end{equation*}

\item[\rm{(\textbf{iv})}] if $A+D \leq \pi$, $B \geq \frac{\pi}{2}$, $C \geq \frac{\pi}{2}$, then $t^2 \geq 0$ and
\begin{equation*}
\mathrm{Re}\, \left(\sinh^{-1}\frac{a_{+}}{t}-\sinh^{-1}\frac{b}{t}-\sinh^{-1}\frac{c}{t}-\sinh^{-1}\frac{a_{-}}{t}\right) = 0,
\end{equation*}

\item[\rm{(\textbf{v})}] if $A+D \leq \pi$, $B \leq \frac{\pi}{2}$, $C \geq \frac{\pi}{2}$ and $t^2 \leq 0$, then
\begin{equation*}
\mathrm{Re}\, \left(\sinh^{-1}\frac{a_{+}}{t}+\sinh^{-1}\frac{b}{t}-\sinh^{-1}\frac{c}{t}-\sinh^{-1}\frac{a_{-}}{t}\right) = 0,
\end{equation*}

\item[\rm{(\textbf{vi})}] if $A+D \leq \pi$, $B \leq \frac{\pi}{2}$, $C \leq \frac{\pi}{2}$, then $t^2 \geq 0$ and
\begin{equation*}
\mathrm{Re}\, \left(\sinh^{-1}\frac{a_{+}}{t}+\sinh^{-1}\frac{b}{t}+\sinh^{-1}\frac{c}{t}-\sinh^{-1}\frac{a_{-}}{t}\right) = 0,
\end{equation*}
\end{itemize}
\end{proposition}
\begin{proof}
Consider case \rm{(\textbf{i})}. For the edge matrix $G^{\star}_s$ of the associated symmetric tetrahedron the following equality holds:
\begin{equation*}
\frac{c}{a_{-}} s^{\star}_{00} + \frac{b}{a_{-}} s^{\star}_{01} + \frac{a_{+}}{a_{-}} s^{\star}_{02} + s^{\star}_{03} = 0.
\end{equation*}
Proposition~\ref{Inequality_Gram} and Lemma~\ref{lemma_symmetric} imply that
\begin{equation*}
s^{\star}_{00} \geq 0,\, s^{\star}_{01}\,\mathcal{A}_{+} \leq 0,\, s^{\star}_{02}\,\mathcal{B} \leq 0,\, s^{\star}_{03}\,\mathcal{C} \leq 0,
\end{equation*}
where all the quantities
\begin{equation*}
\mathcal{A}_{+} = \cos \frac{A+D}{2},\, \mathcal{B} = \cos B,\, \mathcal{C} = \cos C
\end{equation*}
are non-positive under assumptions of case \rm{(\textbf{i})}. Meanwhile $\mathcal{A}_{-} = \cos \frac{D-A}{2}$ is non-negative.
Therefore,
\begin{equation*}
s^{\star}_{00} \geq 0,\, s^{\star}_{01} \geq 0,\, s^{\star}_{02} \geq 0,\, s^{\star}_{03} \geq 0
\end{equation*}
and, by Lemma~\ref{lemma_square_roots},
\begin{equation*}
\sqrt{a^2_{-} - t^2} = a^3_{-} \frac{s^{\star}_{00}}{\sqrt{\Delta}},\, \sqrt{a^2_{+} - t^2} = a^3_{-} \frac{s^{\star}_{01}}{\sqrt{\Delta}},
\end{equation*}
\begin{equation*}
\sqrt{b^2 - t^2} = a^3_{-} \frac{s^{\star}_{02}}{\sqrt{\Delta}},\, \sqrt{c^2 - t^2} = a^3_{-} \frac{s^{\star}_{03}}{\sqrt{\Delta}}.
\end{equation*}

So, the following equality holds:
\begin{equation*}
a_{+} \sqrt{b^2 - t^2} + b \sqrt{a^2_{+} - t^2} = - c \sqrt{a^2_{-} - t^2} - a_{-} \sqrt{c^2 - t^2}.
\end{equation*}
Suppose that $t\neq 0$. Then the equivalent form of the equality above is
\begin{equation*}
\frac{b}{t} \sqrt{\frac{a^2_{+}}{t^2} - 1} + \frac{a_{+}}{t} \sqrt{\frac{b^2}{t^2} - 1} = - \frac{c}{t} \sqrt{\frac{a^2_{-}}{t^2} - 1} - \frac{a_{-}}{t} \sqrt{\frac{c^2}{t^2} - 1}.
\end{equation*}

Applying function $\sinh^{-1}(\circ)$ to the both sides of the equality above and making use of relation \rm{(\textbf{i})} from Lemma~\ref{lemma_inv_hyp_identities}, one obtains equality \rm{(\textbf{i})} of the present Proposition. If $t=0$, then the statement holds in the limiting case $t\rightarrow 0$. The proof for cases \rm{(\textbf{iii})} and \rm{(\textbf{v})} follows by analogy.

Considering cases \rm{(\textbf{ii})}, \rm{(\textbf{iv})} and \rm{(\textbf{vi})}, note that in consequence of the assumptions imposed on parameters $A$, $B$, $C$, $D$ the quantity $\tau$ is purely imaginary and $\mathrm{Im}\,\tau \geq 0$. Then $t$ is also purely imaginary and $\mathrm{Im}\,t \geq 0$. The rest of the proof follows by analogy with cases \rm{(\textbf{i})}, \rm{(\textbf{iii})} and \rm{(\textbf{v})}, making use of Lemma~\ref{lemma_inv_hyp_identities}.
\end{proof}

\begin{proposition}\label{proposition_real_parts_dual}
Let $\mathbf{T}$ be a spherical $\mathbb{Z}_2$-symmetric with dihedral angles $A$, $B = E$, $C = F$, $D$ and edge lengths $l_A$, $l_B = l_E$, $l_C = l_F$, $l_D$. Without loss of generality, assume that $A \geq D$ or, equivalently, $l_A \leq l_D$ and, furthermore, $l_B \geq l_C$.

Then in cases
\begin{itemize}
\item[$\mathrm{(\mathbf{i})}^{\star}$] $A+D \geq \pi$, $B \geq \frac{\pi}{2}$, $C \geq \frac{\pi}{2}$ and $t^2 \geq 0$,
\item[$\mathrm{(\mathbf{ii})}^{\star}$] $A+D \geq \pi$, $B \leq \frac{\pi}{2}$, $C \leq \frac{\pi}{2}$ and $t^2 \geq 0$,
\item[$\mathrm{(\mathbf{iii})}^{\star}$] $A+D \leq \pi$, $B \leq \frac{\pi}{2}$, $C \geq \frac{\pi}{2}$ and $t^2 \geq 0$,
\end{itemize}
Proposition \ref{proposition_real_parts} holds for the tetrahedron $\mathbf{T}^{\star}$ that is dual to the given one.
\end{proposition}
\begin{proof}
By means of the equality
\begin{equation*}
t^2 = \frac{4\,(a_{+}a_{-}-bc)(a_{+}b-a_{-}c)(a_{+}c-a_{-}b)}{\Delta},
\end{equation*}
with $a_{+} = \cos \frac{l_A+l_D}{2}$, $b = \cos l_B$, $c = \cos l_C$, $a_{-} = \cos \frac{l_D-l_A}{2}$, the parameter $t$ can be real only if not all of the quantities $a_{+}$, $b$, $c$ are negative.

Without loss of generality, assume that the following cases are possible:
\begin{itemize}
\item[$\mathrm{(\mathbf{i})}^{\star\star}$] $l_A+l_D \leq \pi$, $l_B \leq \frac{\pi}{2}$, $l_C \leq \frac{\pi}{2}$,
\item[$\mathrm{(\mathbf{ii})}^{\star\star}$] $l_A+l_D \leq \pi$, $l_B \geq \frac{\pi}{2}$, $l_C \geq \frac{\pi}{2}$,
\item[$\mathrm{(\mathbf{iii})}^{\star\star}$] $l_A+l_D \geq \pi$, $l_B \geq \frac{\pi}{2}$, $l_C \leq \frac{\pi}{2}$.
\end{itemize}
Each case above implies that the dihedral angles of the dual tetrahedron $\mathbf{T}^{\star}$ fall under conditions \textrm{(\textbf{i})}, \textrm{(\textbf{iii})} or \textrm{(\textbf{v})} of Proposition \ref{proposition_real_parts}. Parameter $\tau^{\star}$ of the tetrahedron $\mathbf{T}^{\star}$ computed from its dihedral angled satisfies the equality $(\tau^{\star})^2 = -t^2 \leq 0$. It implies that the parameter $t^{\star}$ for the dual tetrahedron $\mathbf{T}^{\star}$ computed from its edge lengths also satisfies condition $(t^{\star})^2 \leq 0$.

Thus, the tetrahedron $\mathbf{T}^{\star}$, that is dual to the given one, falls under one of the cases \textrm{(\textbf{i})}, \textrm{(\textbf{iii})}, \textrm{(\textbf{v})} of Proposition \ref{proposition_real_parts}.
\end{proof}

\subsection{Volume of a $\mathbb{Z}_2$-symmetric spherical tetrahedron}

Let $\mathbf{T}$ be a $\mathbb{Z}_2$-symmetric spherical tetrahedron with dihedral angles $A$, $B = E$, $C = F$, $D$ and edge lengths $l_A$, $l_B = l_E$, $l_C = l_F$, $l_D$. Denote
\begin{equation*}
l_A^+ = \frac{l_A+l_D}{2},\, l_A^- = \frac{l_A-l_D}{2},\, A_{+} = \frac{A+D}{2},\, A_{-} = \frac{D-A}{2},
\end{equation*}
\begin{equation*}
a_{+}=\cos l_A^+,\, a_{-}=\cos l_A^-,\, b=\cos l_B,\, c=\cos l_C.
\end{equation*}
Recall that the principal parameter $u$ of the tetrahedron $\mathbf{T}$ is the positive root of quadratic equation
\begin{equation*}
u^2 + \frac{4\,(a_{+}a_{-}-bc)(a_{+}b-a_{-}c)(a_{+}c-a_{-}b)}{\Delta^{\star}} = 1,
\end{equation*}
with
\begin{equation*}
\Delta^{\star} = (a_{+}+a_{-}+b+c)(a_{+}+a_{-}-b-c)(a_{+}-a_{-}-b+c)(a_{+}-a_{-}+b-c).
\end{equation*}
The auxiliary parameter $t$ from Proposition \ref{proposition_real_parts} satisfies the equality
\begin{equation*}
t^2 = 1 - u^2 = \frac{4\,(a_{+}a_{-}-bc)(a_{+}b-a_{-}c)(a_{+}c-a_{-}b)}{\Delta^{\star}}.
\end{equation*}
Without loss of generality, distinguish the following cases:
\begin{itemize}
\item[\textrm{(\textbf{i})}] $A_{+} \geq \frac{\pi}{2}$, $B \geq \frac{\pi}{2}$, $C \geq \frac{\pi}{2}$ and $t^2 \leq 0$,
\item[$\mathrm{(\mathbf{i})^{\star}}$] $A_{+} \geq \frac{\pi}{2}$, $B \geq \frac{\pi}{2}$, $C \geq \frac{\pi}{2}$ and $t^2 \geq 0$,
\item[\textrm{(\textbf{ii})}] $A_{+} \geq \frac{\pi}{2}$, $B \leq \frac{\pi}{2}$, $C \geq \frac{\pi}{2}$,
\item[\textrm{(\textbf{iii})}] $A_{+} \geq \frac{\pi}{2}$, $B \leq \frac{\pi}{2}$, $C \leq \frac{\pi}{2}$ and $t^2 \leq 0$,
\item[$\mathrm{(\mathbf{iii})^{\star}}$] $A_{+} \geq \frac{\pi}{2}$, $B \leq \frac{\pi}{2}$, $C \leq \frac{\pi}{2}$ and $t^2 \geq 0$,
\item[\textrm{(\textbf{iv})}] $A_{+} \leq \frac{\pi}{2}$, $B \geq \frac{\pi}{2}$, $C \geq \frac{\pi}{2}$,
\item[\rm{(\textbf{v})}]  $A_{+} \leq \frac{\pi}{2}$, $B \leq \frac{\pi}{2}$, $C \geq \frac{\pi}{2}$ and $t^2 \leq 0$,
\item[$\mathrm{(\mathbf{v})^{\star}}$]  $A_{+} \leq \frac{\pi}{2}$, $B \leq \frac{\pi}{2}$, $C \geq \frac{\pi}{2}$ and $t^2 \geq 0$,
\item[\textrm{(\textbf{vi})}] $A_{+} \leq \frac{\pi}{2}$, $B \leq \frac{\pi}{2}$, $C \leq \frac{\pi}{2}$.
\end{itemize}

Define the auxiliary function
\begin{equation*}
\mathrm{V}(\ell, u)  = \frac{1}{2} \int_{\ell}^{\pi/2} \mathrm{Im}\, \log \frac{1 - i \sqrt{u^2/\sin^2 \sigma-1}}{1 + i \sqrt{u^2/\sin^2 \sigma-1}} \, \mathrm{d} \sigma.
\end{equation*}
for all $(\ell, u) \in \mathbb{R}^2$. The branch cut of $\log(\circ)$ runs from $-\infty$ to $0$. The detailed properties of the the function $\mathrm{V}(\circ,\, \circ)$ will be specified in the next section.

Set 
\begin{equation*}
\mathcal{H} = \left(\frac{\pi}{2} - A_+\right)l_A^+ + \left(\frac{\pi}{2} - B\right)l_B + \left(\frac{\pi}{2} - C\right)l_C - \left(\frac{\pi}{2} - A_-\right)l_A^-
\end{equation*}
and
\begin{equation*}
\mathcal{I} = \mathrm{sgn}\left(\frac{\pi}{2} - A_+\right)\mathrm{V}(l_A^+, u) + \mathrm{sgn}\left(\frac{\pi}{2} - B\right)\mathrm{V}(l_B, u) 
\end{equation*}
\begin{equation*}
+ \mathrm{sgn}\left(\frac{\pi}{2} - C\right)\mathrm{V}(l_C, u) - \mathrm{V}(l_A^-, u),
\end{equation*}
where $\mathrm{sgn}(\circ)$ means the sign function.

The following theorem takes place:
\begin{theorem}\label{thm_tetr}
Let $\mathbf{T}$ be a $\mathbb{Z}_2$-symmetric spherical tetrahedron with dihedral angles $A$, $B = E$, $C = F$, $D$ and edge lengths $l_A$, $l_B = l_E$, $l_C = l_F$, $l_D$. Without loss of generality, assume that $A \leq D$ or, equivalently, $l_A \geq l_D$ and, furthermore, $B \leq C$. Then in case tetrahedron $\mathbf{T}$ satisfies the condition $t^2 \leq 0$ its volume is given by the formula
\begin{equation*}
\mathrm{Vol}\,\mathbf{T} = \mathcal{I} - \mathcal{H}.
\end{equation*}
\end{theorem}
\begin{proof}
To prove the theorem we need to show that
\begin{itemize}
\item[\rm{(\textbf{i})}] the function $\mathrm{Vol}\,\mathbf{T}$ satisfies the Schl\"afli formula from Theorem \ref{Schlaefli_formula},
\item[\rm{(\textbf{ii})}] the function $\mathrm{Vol}\,\mathbf{T}$ for the tetrahedron $\mathbf{T}$ with edge lengths $l_A=l_B=l_C=l_D=\pi/2$ equals $\pi^2/8$.
\end{itemize}
Subject to the condition of the Theorem, there are possible cases \rm{(\textbf{i})}-\rm{(\textbf{vi})} pointed above. Consider case \rm{(\textbf{i})}:  $\pi/2 \leq A_+ \leq \pi$, $\pi/2 \leq B \leq \pi$, $\pi/2 \leq C \leq \pi$. By the assumption of the Theorem one has $0\leq A_- \leq \pi/2$. Thus,
\begin{equation*}
\mathcal{I} = -\mathrm{V}(l_A^+, u) - \mathrm{V}(l_B, u) - \mathrm{V}(l_C, u) - \mathrm{V}(l_A^-, u)
\end{equation*}
and
\begin{equation*}
\mathcal{H} = \left(\frac{\pi}{2} - A_+\right)l_A^+ + \left(\frac{\pi}{2} - B\right)l_B + \left(\frac{\pi}{2} - C\right)l_C - \left(\frac{\pi}{2} - A_-\right)l_A^-.
\end{equation*}

Note that if $u \geq 0$ and $0 \leq \ell \leq \pi$ then
\begin{equation*}
\mathrm{V}(\ell, u) = \int_\ell^{\pi/2} \mathrm{Re} \sin^{-1} \frac{\sin \sigma}{u}\, \mathrm{d}\sigma + \frac{\pi}{2}\left(\ell-\frac{\pi}{2}\right),
\end{equation*}
where the branch cut of $\sin^{-1}(\circ)$ is $(-\infty, -1)\cup(1, \infty)$.

It follows that the considered function equals
\begin{equation*}
\mathrm{Vol}\,\mathbf{T} = \mathrm{I} + \mathrm{H} + \pi^2,
\end{equation*}
where
\begin{equation*}
\mathrm{I} = -\int_{l_A^+}^{\pi/2} \mathrm{Re} \sin^{-1} \frac{\sin \sigma}{u}\, \mathrm{d}\sigma -\int_{l_B}^{\pi/2} \mathrm{Re} \sin^{-1} \frac{\sin \sigma}{u}\, \mathrm{d}\sigma 
\end{equation*}
\begin{equation*}
-\int_{l_C}^{\pi/2} \mathrm{Re} \sin^{-1} \frac{\sin \sigma}{u}\, \mathrm{d}\sigma - \int_{l_A^-}^{\pi/2} \mathrm{Re} \sin^{-1} \frac{\sin \sigma}{u}\, \mathrm{d}\sigma - \pi l_A^+ - \pi l_B - \pi l_C,
\end{equation*}
and
\begin{equation*}
\mathrm{H} = A^+ l_A^+ + B l_B + C l_C - A^- l_A^-
\end{equation*}
\begin{equation*}
= \frac{1}{2} A l_A + B l_B + C l_C + \frac{1}{2} D l_D.
\end{equation*}

Once we prove
\begin{equation*}
\mathrm{d}\mathrm{I} = - \frac{1}{2} A \mathrm{d}l_A - B \mathrm{d}l_B - C \mathrm{d}l_C - \frac{1}{2} D \mathrm{d}l_D
\end{equation*}
it follows that
\begin{equation*}
\mathrm{d}\mathrm{Vol}\,\mathbf{T} = \frac{1}{2} l_A \mathrm{d}A + l_B \mathrm{d}B + l_C \mathrm{d}C + \frac{1}{2} l_D \mathrm{d}D
\end{equation*}
and condition \rm{(\textbf{i})} is fulfilled.

Compute the partial derivative 
\begin{equation*}
\frac{\partial \mathrm{I}}{\partial l_A} = -\frac{1}{2} \sin^{-1} \frac{\sin l_A^+}{u} + \frac{1}{2} \sin^{-1} \frac{\sin l_A^-}{u} - \frac{\pi}{2} + \frac{1}{u} \frac{\partial u}{\partial l_A} \mathrm{F}(l_A^+,l_B,l_C,l_A^-,u)
\end{equation*}
where
\begin{equation*}
\mathrm{F}(l_A^+,l_B,l_C,l_A^-,u) = \mathrm{Re}\, \left(\sinh^{-1}\frac{\cos l_A^+}{\sqrt{1-u^2}}+\sinh^{-1}\frac{\cos l_B}{\sqrt{1-u^2}}\right.
\end{equation*}
\begin{equation*}
\left.+\sinh^{-1}\frac{\cos l_C}{\sqrt{1-u^2}}+\sinh^{-1}\frac{\cos l_A^-}{\sqrt{1-u^2}}\right).
\end{equation*}
Proposition \ref{proposition_real_parts} implies that $\mathrm{F}(l_A^+,l_B,l_C,l_A^-,u) = 0$. Then, by Proposition \ref{SineRule}, the following equality hold:
\begin{equation*}
\frac{\partial \mathrm{I}}{\partial l_A} = -\frac{1}{2} \sin^{-1} \frac{\sin l_A^+}{u} + \frac{1}{2} \sin^{-1} \frac{\sin l_A^-}{u} - \frac{\pi}{2} = 
\end{equation*}
\begin{equation*}
= \frac{1}{2} \sin^{-1} \sin A_+ + \frac{1}{2} \sin^{-1} \sin A_- - \frac{\pi}{2} = 
\end{equation*}
\begin{equation*}
= \frac{1}{2}\left(\pi - \frac{A+D}{2}\right) + \frac{1}{2} \frac{D-A}{2} - \frac{\pi}{2} = -\frac{A}{2},
\end{equation*}
taking into account that
\begin{equation*}
\sin^{-1} x = \left\{ \begin{array}{cc} x, &\mbox{if } 0\leq x\leq \pi/2,\\ \pi - x, &\mbox{if } \pi/2 \leq x\leq\pi. \end{array} \right.
\end{equation*}
Analogously,
\begin{equation*}
\frac{\partial \mathrm{I}}{\partial l_B} = -B,\,\frac{\partial \mathrm{I}}{\partial l_C} = -C,\,\frac{\partial \mathrm{I}}{\partial l_D} = -\frac{D}{2}.
\end{equation*}
Thus, condition \rm{(\textbf{i})} is satisfied.

Compute the function $\mathrm{Vol}\,\mathbf{T}$ with $l_A^{+} = l_B = l_C = \frac{\pi}{2}$, $l_A^{-} = 0$ and $A_{+} = B = C = \frac{\pi}{2}$, $A_{-} = 0$, setting $u = 1$ as follows from Proposition \ref{SineRule}. Then one has $\mathrm{Vol}\,\mathbf{T} = \pi^2/8$ and condition \textrm{(\textbf{ii})} holds. It implies the Theorem for case \textrm{(\textbf{i})} to be proven.

The proof for cases \textrm{(\textbf{ii})}, \textrm{(\textbf{iii})}, \textrm{(\textbf{iv})}, \textrm{(\textbf{v})} and \textrm{(\textbf{vi})} follows by analogy.
\end{proof}

In cases $\mathrm{(\mathbf{i})^{\star}}$, $\mathrm{(\mathbf{iii})^{\star}}$ and $\mathrm{(\mathbf{v})^{\star}}$ the following theorem holds:
\begin{theorem}\label{thm_tetr_dual}
Let $\mathbf{T}$ be a spherical $\mathbb{Z}_2$-symmetric tetrahedron with dihedral angles $A$, $B = E$, $C = F$, $D$ and edge lengths $l_A$, $l_B = l_E$, $l_C = l_F$, $l_D$. Without loss of generality, assume that $A \geq D$ or, equivalently, $l_A \leq l_D$ and, furthermore,  $l_B \geq l_C$. Then, in case tetrahedron $\mathbf{T}$ satisfies the condition $t^2 \geq 0$, the statement of Theorem \ref{thm_tetr} holds for the tetrahedron $\mathbf{T}^{\star}$ that is dual to the given one.
\end{theorem}
\begin{proof}
The proof follows by analogy with Theorem \ref{thm_tetr} using Proposition \ref{proposition_real_parts_dual} instead of Proposition \ref{proposition_real_parts}.
\end{proof}

To find the volume of a tetrahedron $\mathbf{T}$ that respects the conditions of Theorem \ref{thm_tetr_dual} one may apply Theorem \ref{thm_tetr} to the dual tetrahedron $\mathbf{T}^{\star}$ and then make use of the Sforza formula from Theorem \ref{theorem_Sforza}.

\subsection{Computation of certain volumes}

It follows from Lemma \ref{lemma_V_elementary} of the next section that in case $u=1$ the function $\mathrm{V}(\ell, u)$ has rather elementary form. Thus, the volume of a tetrahedron with the principal parameter $u=1$ can be represented by the elementary functions. The equality $u=1$ also means the same as  $t=0$, because of the relation $t^2 = u^2 - 1$ and non-negativity of $u$.

Consider the associated symmetric tetrahedron $\mathbf{T}_s$ with its auxiliary parameter $t_s = 0$ because of the relation $t = a_{-} t_s$ between $t$ and $t_s$ from the proof of Lemma \ref{lemma_symmetric}.

By Lemma \ref{lemma_symmetric} one has
\begin{equation*}
t_s = \frac{4(\tilde{a}-\tilde{b}\tilde{c})(\tilde{b}-\tilde{a}\tilde{c})(\tilde{c}-\tilde{a}\tilde{b})}{\Delta^{\star}},
\end{equation*}
where $\Delta^{\star} = \det G^{\star}_s$ is the determinant of the edge matrix $G^{\star}_s$ of the tetrahedron $\mathbf{T}_s$. Also the following equalities hold:
\begin{equation*}
\tilde{a} = \cos l_\alpha = \frac{a_{+}}{a_{-}},\,\,\tilde{b} = \cos l_\beta = \frac{b}{a_{-}},\,\,\tilde{c} = \cos l_\gamma = \frac{c}{a_{-}}.
\end{equation*}

The equality $t_s = 0$ implies three cases: $\tilde{a} - \tilde{b}\tilde{c}=0$, or $\tilde{b} - \tilde{a}\tilde{c}=0$, or $\tilde{c} - \tilde{a}\tilde{b}=0$. Together, these equalities imply either $\tilde{a} = \tilde{b} = \tilde{c} = \pm 1$ and the tetrahedron is $\mathbf{T}$ degenerate, or $\tilde{a} = \tilde{b} = \tilde{c} = 0$ and both tetrahedra $\mathbf{T}$ and $\mathbf{T}_s$ are isometric to a equilateral tetrahedron with edge length $\frac{\pi}{2}$. 

Without loss of generality, suppose only two of the equalities above to be hold: $\tilde{b} - \tilde{a}\tilde{c} = 0$ and $\tilde{c} - \tilde{a}\tilde{b} = 0$. If the tetrahedron $\mathbf{T}_s$ is not degenerate, then one obtains $\tilde{b} = \tilde{c} = 0$. Therefore, tetrahedra $\mathbf{T}_s$ provide a one-parametric family of tetrahedra with $0 < l_\alpha < \pi$, $l_\beta = l_\gamma = \frac{\pi}{2}$. The associated tetrahedron $\mathbf{T}$ has edge lengths $0 < l_A,\,\,l_D < \pi$, $l_B = l_C = \frac{\pi}{2}$.

Suppose now that only one equality, namely $\tilde{a} - \tilde{b}\tilde{c} = 0$, holds. By Lemma~\ref{lemma_symmetric} and formul{\ae} of spherical geometry from \cite[Ch.~1, \textsection~4.2]{Vinberg} one obtains
\begin{equation*}
\cos \alpha = -\tilde{a},\,\, \cos \beta = \tilde{b},\,\, \cos \gamma = \tilde{c}.
\end{equation*}

Thus, for $\mathbf{T}$ the following inequalities hold:
\begin{equation*}
\cos l_A^{+} \cos A^{+} \leq 0,\,\, \cos l_B \cos B \geq 0,\,\, \cos l_C \cos C \geq 0.
\end{equation*}

Apply Theorem \ref{SineRule} to the tetrahedron $\mathbf{T}$ with principal parameter $u = 1$ and obtain that
\begin{equation*}
\sin l_A^{+} = \sin A^{+},\,\, \sin l_B = \sin B,\,\, \sin l_C = \sin C,\,\, \sin l_A^{-} = \sin A^{-}.
\end{equation*}
From the above one derives the following equalities:
\begin{equation*}
A^{+} = \pi - l_A^{+},\,\, B = l_B,\,\, C = l_C,\,\, A^{-} = l_A^{-}
\end{equation*}
or, equivalently,
\begin{equation*}
A = \pi - l_A,\,\, B = l_B,\,\, C = l_C,\,\, D = \pi - l_D.
\end{equation*}

The cases of equalities $\tilde{b} - \tilde{a}\tilde{c} = 0$ and $\tilde{c} - \tilde{a}\tilde{b} = 0$ are analogous. Moreover, they are the same up to a permutation of the parameters $l_B$ and $l_C$ of the tetrahedron $\mathbf{T}$.

Thus, the other possible equalities are
\begin{equation*}
A = l_D,\,\, B = \pi - l_B,\,\, C = l_C,\,\, D = l_A
\end{equation*}
or
\begin{equation*}
A = l_D,\,\, B = l_B,\,\, C = \pi - l_C,\,\, D = l_A.
\end{equation*}

Note, that the last three considered cases cover all the occasions mentioned above, like equilateral tetrahedron with edge length $\frac{\pi}{2}$ or the family of tetrahedra with edge lengths $0 < l_A,\,\,l_D < \pi$, $l_B = l_C = \frac{\pi}{2}$.

The following statement holds:
\begin{proposition}\label{proposition_elementary}
Let $\mathbf{T}$ be a spherical $\mathbb{Z}_2$-symmetric tetrahedron. Suppose

\begin{equation*}
\cos l_A^+ \cos l_A^- - \cos l_B \cos l_C = 0,
\end{equation*}
or
\begin{equation*}
\cos l_B \cos l_A^- - \cos l_A^+ \cos l_C = 0,
\end{equation*}
or
\begin{equation*}
\cos l_C \cos l_A^- - \cos l_A^+ \cos l_B = 0.
\end{equation*}
Then the volume of $\mathbf{T}$ is given by the corresponding formula
\begin{equation*}
\mathrm{Vol}\,\mathbf{T} = \frac{1}{2}\left(-\frac{l^{2}_A}{2} + l_B^2 + l_C^2 - \frac{l^{2}_D}{2}\right),
\end{equation*}
or
\begin{equation*}
\mathrm{Vol}\,\mathbf{T} = \frac{l_A l_D - l^{2}_B + l_C^2}{2},
\end{equation*}
or
\begin{equation*}
\mathrm{Vol}\,\mathbf{T} = \frac{l_A l_D + l^{2}_B - l_C^2}{2}.
\end{equation*}
\end{proposition}
\begin{proof}
Let us consider the case $\cos l_A^{+} \cos l_A^{-} - \cos l_B \cos l_C = 0$. From the consideration above one obtains that $A = \pi - l_A$, $B = l_B$, $C = l_C$ and $D = \pi - l_D$. The principal parameter $u$ of $\mathbf{T}$ satisfies the equality $u = 1$. Apply Theorem \ref{thm_tetr} and Lemma \ref{lemma_V_elementary} to compute the volume of $\mathbf{T}$ using elementary functions. Simplifying the corresponding equation one arrives at the statement of the Proposition.

The proof for other cases follows by analogy.
\end{proof}

Note, that the claims of Proposition \ref{proposition_elementary} on the tetrahedron $\mathbf{T}$ imply its associated symmetric tetrahedron $\mathbf{T}_s$ has at least one right triangle face. The tetrahedron $\mathbf{T}$ itself might not have such one.

\subsection{Properties of the auxiliary function $\mathrm{V}(\ell, u)$}

The list of basic properties which the function
\begin{equation*}
\mathrm{V}(\ell, u)  = \frac{1}{2} \int_{\ell}^{\pi/2} \mathrm{Im}\, \log \frac{1 - i \sqrt{u^2/\sin^2 \sigma-1}}{1 + i \sqrt{u^2/\sin^2 \sigma-1}} \, \mathrm{d} \sigma
\end{equation*}
with $(\ell, u) \in \mathbb{R}^2$ enjoys is given below:
\begin{lemma}\label{lemma_V_function}
The function $V(\ell, u)$ defined above satisfies the following properties:
\begin{enumerate}
\item[\rm{(\textbf{i})}] $V(\ell, u)$ is continuous and a.e. differentiable in $\mathbb{R}^2$.
\end{enumerate}
For all $(\ell, u) \in \mathbb{R}^2$
\begin{enumerate}
\item[\rm{(\textbf{ii})}] $V(\ell, u) = V(\ell, -u)$,
\item[\rm{(\textbf{iii})}] $V(\pi-\ell, u) = - V(\ell, u)$,
\item[\rm{(\textbf{iv})}] $V(\ell, u) + V(-\ell, u) = 2\,V(0, u)$.
\end{enumerate}
The function $V(\ell, u)$ is linear periodic with respect to $\ell$, that is
\begin{enumerate}
\item[\rm{(\textbf{v})}] $V(\ell + k \pi, u) = V(\ell, u) - 2\,k\,V(0, u)$ for all $k \in \mathbb{Z}$.
\end{enumerate}
\end{lemma}
\begin{proof}
The properties \rm{(\textbf{i})}-\rm{(\textbf{iii})} follow immediately from the definition of the function $\mathrm{V}(\ell, u)$. 

To prove \rm{(\textbf{iv})} notice that the equality holds if $\ell = 0$. In accordance with the definition of $V(\ell, u)$, the derivatives of both sides of \rm{(\textbf{iv})} with respect to $\ell$ vanish. Thus, the equality holds.

The derivatives of both sides of \rm{(\textbf{v})} with respect to $\ell$ are equal. Verification of the equality for $\ell=0$ results in the complete proof of \rm{(\textbf{v})}.
By \rm{(\textbf{iii})} and \rm{(\textbf{iv})} it follows that
\begin{equation*}
V(k\pi+\pi, u) = -V(-k\pi, u),\,\, -V(-k\pi, u) = V(k\pi, u) - 2 V(0, u),
\end{equation*}
with $k \in \mathbb{Z}$.
Hence
\begin{equation*}
V(k\pi+\pi, u) = V(k\pi, u) - 2 V(0, u) = \dots = V(0, u) - 2 (k+1) V(0, u)
\end{equation*}
and equality \rm{(\textbf{v})} holds.
\end{proof}

For the special value of $u = 1$ the function $\mathrm{V}(\ell, 1)$ can be expressed by elementary functions:
\begin{lemma}\label{lemma_V_elementary}
The function $\mathrm{V}(\ell, u)$ with $u = 1$, $0\leq \ell \leq \pi$ can be expressed as
\begin{equation*}
\mathrm{V}(\ell, 1) = \frac{1}{2}\left(\ell-\frac{\pi}{2}\right)\left|\ell-\frac{\pi}{2}\right|.
\end{equation*}
\end{lemma}
\begin{proof}
Notice, that if $u\geq 0$ and $0\leq \ell \leq \pi$ then
\begin{equation*}
V(\ell, u) = \int_{\ell}^{\pi/2}\mathrm{Re}\,\sin^{-1}\frac{\sin\sigma}{u}\,\mathrm{d}\sigma + \frac{\pi}{2}\left( \ell - \frac{\pi}{2} \right).
\end{equation*}
Put $u=1$ and use the equality
\begin{equation*}
\sin^{-1} x = \left\{ \begin{array}{cc} x, &\mbox{if } 0\leq x\leq \pi/2,\\ \pi - x, &\mbox{if } \pi/2 \leq x\leq\pi. \end{array} \right.
\end{equation*}
It follows that
\begin{equation*}
V(\ell, 1) = \left\{ \begin{array}{cc} -1/2 (\ell-\pi/2)^2, &\mbox{if } 0\leq \ell \leq \pi/2,\\ 1/2 (\ell-\pi/2)^2, &\mbox{if } \pi/2 \leq \ell \leq\pi. \end{array} \right.
\end{equation*}
\end{proof}

If $u \geq 1$ then we have the following
\begin{lemma}\label{lemma_V_analytic}
The function $\mathrm{V}(\ell, u)$ with $u \geq 1$ has the series representation
\begin{equation*}
\mathrm{V}(\ell, u) = \frac{\pi}{2}\left(\ell-\frac{\pi}{2}\right) + \sum_{k=0}^\infty p_k\,\frac{ u^{-2k-1}}{(2k+1)^2}
\end{equation*}
with $p_k = 1 - B(\sin \ell; k+1, 1/2)/B(k+1, 1/2)$, where $B(\circ, \circ)$ means beta-function and $B(\circ; \circ, \circ)$ means incomplete beta-function.
\end{lemma}
\begin{proof}
Use the following series representation of the integrand in the expression for $V(u, \phi)$ with respect to the variable $u$ at the point $u = \infty$:
\begin{equation*}
\frac{1}{2}\,\mathrm{Im}\, \log \frac{1 - i \sqrt{u^2/\sin^2 \sigma-1}}{1 + i \sqrt{u^2/\sin^2 \sigma-1}} = \frac{\pi}{2} + \sum_{k=0}^{\infty} \frac{(2 k + 1)!!}{k!\,2^k\,(2k+1)^2} \left( \frac{\sin \sigma}{u} \right)^{2k+1}.
\end{equation*}
The representation above holds for all $u\in [1,\infty)$ and $\sigma \in \mathbb{R}$. Integrating the series above with respect to $\sigma$ from $\phi$ to $\frac{\pi}{2}$ with $\phi \in [0, \pi]$ finishes the proof.
\end{proof}

If $u \leq 1$ then the function $\mathrm{V}(\ell, u)$ has discontinuous second partial derivatives. Their points of discontinuity in the set $\{(\phi, u)\in[0,\pi]\times (0,1)\}$ are $\left(\frac{\pi}{2}\pm(\frac{\pi}{2}-\sin^{-1}u),\,u\right)$.

\flushleft{\emph{
Alexander Kolpakov \\
Universit\"at Freiburg\\
Fachbereich Mathematik\\
Museengasse, Geb. 23,\\
Freiburg, CH-1700 Schweiz\\
}
\rm{aleksandr.kolpakov@unifr.ch}

\medskip
\flushleft{\emph{
Alexander Mednykh\\
Sobolev Institute of Mathematics, SB RAS\\
630090, Koptyug avenue, bld. 4,\\
Novosibirsk, Russia\\
}
\rm{mednykh@math.nsc.ru}
}

\medskip
\flushleft{\emph{
Marina Pa\v{s}kevi\v{c}\\
Novosibirsk State University\\
630090, Pirogova str., bld. 2,\\
Novosibirsk, Russia\\
}
\rm{Pashkevich\_M@mail.ru}
}

\end{document}